\expandafter\chardef\csname pre amssym.def at\endcsname=\the\catcode`\@
\catcode`\@=11

\def\undefine#1{\let#1\undefined}
\def\newsymbol#1#2#3#4#5{\let\next@\relax
 \ifnum#2=\@ne\let\next@\msafam@\else
 \ifnum#2=\tw@\let\next@\msbfam@\fi\fi
 \mathchardef#1="#3\next@#4#5}
\def\mathhexbox@#1#2#3{\relax
 \ifmmode\mathpalette{}{\m@th\mathchar"#1#2#3}%
 \else\leavevmode\hbox{$\m@th\mathchar"#1#2#3$}\fi}
\def\hexnumber@#1{\ifcase#1 0\or 1\or 2\or 3\or 4\or 5\or 6\or 7\or 8\or
 9\or A\or B\or C\or D\or E\or F\fi}

\font\tenmsa=msam10
\font\sevenmsa=msam7
\font\fivemsa=msam5
\newfam\msafam
\textfont\msafam=\tenmsa
\scriptfont\msafam=\sevenmsa
\scriptscriptfont\msafam=\fivemsa
\edef\msafam@{\hexnumber@\msafam}
\mathchardef\dabar@"0\msafam@39
\def\dashrightarrow{\mathrel{\dabar@\dabar@\mathchar"0\msafam@4B}}
\def\dashleftarrow{\mathrel{\mathchar"0\msafam@4C\dabar@\dabar@}}

\def\ulcorner{\delimiter"4\msafam@70\msafam@70 }
\def\urcorner{\delimiter"5\msafam@71\msafam@71 }
\def\llcorner{\delimiter"4\msafam@78\msafam@78 }
\def\lrcorner{\delimiter"5\msafam@79\msafam@79 }
\def\yen{{\mathhexbox@\msafam@55 }}
\def\checkmark{{\mathhexbox@\msafam@58 }}
\def\circledR{{\mathhexbox@\msafam@72 }}
\def\maltese{{\mathhexbox@\msafam@7A }}

\font\tenmsb=msbm10
\font\sevenmsb=msbm7
\font\fivemsb=msbm5
\newfam\msbfam
\textfont\msbfam=\tenmsb
\scriptfont\msbfam=\sevenmsb
\scriptscriptfont\msbfam=\fivemsb
\edef\msbfam@{\hexnumber@\msbfam}

\catcode`\@=\csname pre amssym.def at\endcsname

\expandafter\ifx\csname pre amssym.tex at\endcsname\relax \else \endinput\fi
\expandafter\chardef\csname pre amssym.tex at\endcsname=\the\catcode`\@
\catcode`\@=11
\newsymbol\boxdot 1200
\newsymbol\boxplus 1201
\newsymbol\boxtimes 1202
\newsymbol\square 1003
\newsymbol\blacksquare 1004
\newsymbol\centerdot 1205
\newsymbol\lozenge 1006
\newsymbol\blacklozenge 1007
\newsymbol\circlearrowright 1308
\newsymbol\circlearrowleft 1309
\undefine\rightleftharpoons
\newsymbol\rightleftharpoons 130A
\newsymbol\leftrightharpoons 130B
\newsymbol\boxminus 120C
\newsymbol\Vdash 130D
\newsymbol\Vvdash 130E
\newsymbol\vDash 130F
\newsymbol\twoheadrightarrow 1310
\newsymbol\twoheadleftarrow 1311
\newsymbol\leftleftarrows 1312
\newsymbol\rightrightarrows 1313
\newsymbol\upuparrows 1314
\newsymbol\downdownarrows 1315
\newsymbol\upharpoonright 1316
 
\newsymbol\downharpoonright 1317
\newsymbol\upharpoonleft 1318
\newsymbol\downharpoonleft 1319
\newsymbol\rightarrowtail 131A
\newsymbol\leftarrowtail 131B
\newsymbol\leftrightarrows 131C
\newsymbol\rightleftarrows 131D
\newsymbol\Lsh 131E
\newsymbol\Rsh 131F
\newsymbol\rightsquigarrow 1320
\newsymbol\leftrightsquigarrow 1321
\newsymbol\looparrowleft 1322
\newsymbol\looparrowright 1323
\newsymbol\circeq 1324
\newsymbol\succsim 1325
\newsymbol\gtrsim 1326
\newsymbol\gtrapprox 1327
\newsymbol\multimap 1328
\newsymbol\therefore 1329
\newsymbol\because 132A
\newsymbol\doteqdot 132B
 
\newsymbol\triangleq 132C
\newsymbol\precsim 132D
\newsymbol\lesssim 132E
\newsymbol\lessapprox 132F
\newsymbol\eqslantless 1330
\newsymbol\eqslantgtr 1331
\newsymbol\curlyeqprec 1332
\newsymbol\curlyeqsucc 1333
\newsymbol\preccurlyeq 1334
\newsymbol\leqq 1335
\newsymbol\leqslant 1336
\newsymbol\lessgtr 1337
\newsymbol\backprime 1038
\newsymbol\risingdotseq 133A
\newsymbol\fallingdotseq 133B
\newsymbol\succcurlyeq 133C
\newsymbol\geqq 133D
\newsymbol\geqslant 133E
\newsymbol\gtrless 133F
\newsymbol\sqsubset 1340
\newsymbol\sqsupset 1341
\newsymbol\vartriangleright 1342
\newsymbol\vartriangleleft 1343
\newsymbol\trianglerighteq 1344
\newsymbol\trianglelefteq 1345
\newsymbol\bigstar 1046
\newsymbol\between 1347
\newsymbol\blacktriangledown 1048
\newsymbol\blacktriangleright 1349
\newsymbol\blacktriangleleft 134A
\newsymbol\vartriangle 134D
\newsymbol\blacktriangle 104E
\newsymbol\triangledown 104F
\newsymbol\eqcirc 1350
\newsymbol\lesseqgtr 1351
\newsymbol\gtreqless 1352
\newsymbol\lesseqqgtr 1353
\newsymbol\gtreqqless 1354
\newsymbol\Rrightarrow 1356
\newsymbol\Lleftarrow 1357
\newsymbol\veebar 1259
\newsymbol\barwedge 125A
\newsymbol\doublebarwedge 125B
\undefine\angle
\newsymbol\angle 105C
\newsymbol\measuredangle 105D
\newsymbol\sphericalangle 105E
\newsymbol\varpropto 135F
\newsymbol\smallsmile 1360
\newsymbol\smallfrown 1361
\newsymbol\Subset 1362
\newsymbol\Supset 1363
\newsymbol\Cup 1264
 
\newsymbol\Cap 1265
 
\newsymbol\curlywedge 1266
\newsymbol\curlyvee 1267
\newsymbol\leftthreetimes 1268
\newsymbol\rightthreetimes 1269
\newsymbol\subseteqq 136A
\newsymbol\supseteqq 136B
\newsymbol\bumpeq 136C
\newsymbol\Bumpeq 136D
\newsymbol\lll 136E
 
\newsymbol\ggg 136F
 
\newsymbol\circledS 1073
\newsymbol\pitchfork 1374
\newsymbol\dotplus 1275
\newsymbol\backsim 1376
\newsymbol\backsimeq 1377
\newsymbol\complement 107B
\newsymbol\intercal 127C
\newsymbol\circledcirc 127D
\newsymbol\circledast 127E
\newsymbol\circleddash 127F
\newsymbol\lvertneqq 2300
\newsymbol\gvertneqq 2301
\newsymbol\nleq 2302
\newsymbol\ngeq 2303
\newsymbol\nless 2304
\newsymbol\ngtr 2305
\newsymbol\nprec 2306
\newsymbol\nsucc 2307
\newsymbol\lneqq 2308
\newsymbol\gneqq 2309
\newsymbol\nleqslant 230A
\newsymbol\ngeqslant 230B
\newsymbol\lneq 230C
\newsymbol\gneq 230D
\newsymbol\npreceq 230E
\newsymbol\nsucceq 230F
\newsymbol\precnsim 2310
\newsymbol\succnsim 2311
\newsymbol\lnsim 2312
\newsymbol\gnsim 2313
\newsymbol\nleqq 2314
\newsymbol\ngeqq 2315
\newsymbol\precneqq 2316
\newsymbol\succneqq 2317
\newsymbol\precnapprox 2318
\newsymbol\succnapprox 2319
\newsymbol\lnapprox 231A
\newsymbol\gnapprox 231B
\newsymbol\nsim 231C
\newsymbol\ncong 231D
\newsymbol\diagup 231E
\newsymbol\diagdown 231F
\newsymbol\varsubsetneq 2320
\newsymbol\varsupsetneq 2321
\newsymbol\nsubseteqq 2322
\newsymbol\nsupseteqq 2323
\newsymbol\subsetneqq 2324
\newsymbol\supsetneqq 2325
\newsymbol\varsubsetneqq 2326
\newsymbol\varsupsetneqq 2327
\newsymbol\subsetneq 2328
\newsymbol\supsetneq 2329
\newsymbol\nsubseteq 232A
\newsymbol\nsupseteq 232B
\newsymbol\nparallel 232C
\newsymbol\nmid 232D
\newsymbol\nshortmid 232E
\newsymbol\nshortparallel 232F
\newsymbol\nvdash 2330
\newsymbol\nVdash 2331
\newsymbol\nvDash 2332
\newsymbol\nVDash 2333
\newsymbol\ntrianglerighteq 2334
\newsymbol\ntrianglelefteq 2335
\newsymbol\ntriangleleft 2336
\newsymbol\ntriangleright 2337
\newsymbol\nleftarrow 2338
\newsymbol\nrightarrow 2339
\newsymbol\nLeftarrow 233A
\newsymbol\nRightarrow 233B
\newsymbol\nLeftrightarrow 233C
\newsymbol\nleftrightarrow 233D
\newsymbol\divideontimes 223E
\newsymbol\varnothing 203F
\newsymbol\nexists 2040
\newsymbol\Finv 2060
\newsymbol\Game 2061
\newsymbol\mho 2066
\newsymbol\eth 2067
\newsymbol\eqsim 2368
\newsymbol\beth 2069
\newsymbol\gimel 206A
\newsymbol\daleth 206B
\newsymbol\lessdot 236C
\newsymbol\gtrdot 236D
\newsymbol\ltimes 226E
\newsymbol\rtimes 226F
\newsymbol\shortmid 2370
\newsymbol\shortparallel 2371
\newsymbol\smallsetminus 2272
\newsymbol\thicksim 2373
\newsymbol\thickapprox 2374
\newsymbol\approxeq 2375
\newsymbol\succapprox 2376
\newsymbol\precapprox 2377
\newsymbol\curvearrowleft 2378
\newsymbol\curvearrowright 2379
\newsymbol\digamma 207A
\newsymbol\varkappa 207B
\newsymbol\Bbbk 207C
\newsymbol\hslash 207D
\undefine\hbar
\newsymbol\hbar 207E
\newsymbol\backepsilon 237F
\catcode`\@=\csname pre amssym.tex at\endcsname

\magnification=1200
\hsize=468truept
\vsize=646truept
\voffset=-10pt
\parskip=4pt
\baselineskip=14truept
\count0=1

\dimen100=\hsize

\def\leftill#1#2#3#4{
\medskip
\line{$
\vcenter{
\hsize = #1truept \hrule\hbox{\vrule\hbox to  \hsize{\hss \vbox{\vskip#2truept
\hbox{{\copy100 \the\count105}: #3}\vskip2truept}\hss }
\vrule}\hrule}
\dimen110=\dimen100
\advance\dimen110 by -36truept
\advance\dimen110 by -#1truept
\hss \vcenter{\hsize = \dimen110
\medskip
\noindent { #4\par\medskip}}$}
\advance\count105 by 1
}
\def\rightill#1#2#3#4{
\medskip
\line{
\dimen110=\dimen100
\advance\dimen110 by -36truept
\advance\dimen110 by -#1truept
$\vcenter{\hsize = \dimen110
\medskip
\noindent { #4\par\medskip}}
\hss \vcenter{
\hsize = #1truept \hrule\hbox{\vrule\hbox to  \hsize{\hss \vbox{\vskip#2truept
\hbox{{\copy100 \the\count105}: #3}\vskip2truept}\hss }
\vrule}\hrule}
$}
\advance\count105 by 1
}
\def\midill#1#2#3{\medskip
\line{$\hss
\vcenter{
\hsize = #1truept \hrule\hbox{\vrule\hbox to  \hsize{\hss \vbox{\vskip#2truept
\hbox{{\copy100 \the\count105}: #3}\vskip2truept}\hss }
\vrule}\hrule}
\dimen110=\dimen100
\advance\dimen110 by -36truept
\advance\dimen110 by -#1truept
\hss $}
\advance\count105 by 1
}
\def\insectnum{\copy110\the\count120
\advance\count120 by 1
}

\font\ninerm=cmr9
\font\eightrm=cmr8

\font\tenrm=cmr10 at 10pt

\font\sc=cmcsc10

\def\msb{\fam\msbfam\tenmsb}

\def\bbc{{\msb C}}

\def\bbi{{\msb I}}

\def\bbp{{\msb P}}
\def\bbq{{\msb Q}}
\def\bbr{{\msb R}}

\def\bbz{{\msb Z}}

\def\grD{\Delta}

\def\grG{\Gamma}

\def\grL{\Lambda}

\def\grS{\Sigma}

\def\gra{\alpha}

\def\gri{\iota}

\def\grl{\lambda}

\def\gro{\omega}

\def\grr{\rho}

\def\grt{\tau}

\def\grz{\zeta}

\def\la#1{\hbox to #1pc{\leftarrowfill}}
\def\ra#1{\hbox to #1pc{\rightarrowfill}}

\def\fract#1#2{\raise4pt\hbox{$ #1 \atop #2 $}}
\def\decdnar#1{\phantom{\hbox{$\scriptstyle{#1}$}}
\left\downarrow\vbox{\vskip15pt\hbox{$\scriptstyle{#1}$}}\right.}

\def\bowtie{\hbox to 1pt{\hss}\raise.66pt\hbox{$\scriptstyle{>}$}
\kern-4.9pt\triangleleft}
\def\hsmash{\triangleright\kern-4.4pt\raise.66pt\hbox{$\scriptstyle{<}$}}
\def\boxit#1{\vbox{\hrule\hbox{\vrule\kern3pt
\vbox{\kern3pt#1\kern3pt}\kern3pt\vrule}\hrule}}

\def\za{\vrule height6pt width4pt depth1pt}

\font\aa=eufm10

\def\Got#1{\hbox{\aa#1}}

\def\bfw{{\bf w}}

\def\calo{{\cal O}}

\def\cald{{\cal D}}

\def\call{{\cal L}}

\def\calo{{\cal O}}

\def\calz{{\cal Z}}

\def\gA{{\Got A}}

\def\gG{{\Got G}}

\def\gI{{\Got I}}


\font\svtnrm=cmr17

\font\bsc=cmcsc10 at 10truept

\def\coker{\hbox{coker}}
\def\div{\hbox{div}~}
\def\lcm{\hbox{lcm}}
\def\Ric{\hbox{Ric}}
\def\Se{Sasakian-Einstein }
\def\dim{\hbox{dim}}
\def\linkint{0}
\def\linksas{1}
\def\linkke{2}
\def\linkex{3}
\def\linkmi{4}
\def\linktop{5}
\def\linkpr{6}

\centerline{\svtnrm New Einstein Metrics in Dimension Five}

\bigskip
\centerline{\sc Charles P. Boyer~~ Krzysztof Galicki~~}
\bigskip
\footnote{}{\ninerm During the preparation of this work the authors 
were partially supported by NSF grant DMS-9970904.}
\bigskip

\centerline{\vbox{\hsize = 5.85truein
\baselineskip = 12.5truept
\eightrm
\noindent {\bsc Abstract:}
The purpose of this note is to introduce a new method for proving the existence of \Se 
metrics on certain simply connected odd dimensional manifolds. We then apply this 
method to prove the existence of new \Se
metrics on $\scriptstyle{S^2\times S^3}$ and on 
$\scriptstyle{(S^2\times S^3)\# (S^2\times S^3).}$ These
give the first known examples of non-regular \Se 5-manifolds. Our method
involves describing the \Se structures as links of certain isolated
hypersurface singularities, and makes use of the recent work of Demailly and
Koll\'ar who  obtained new examples of K\"ahler-Einstein del Pezzo surfaces
with quotient singularities. }} \tenrm

\bigskip
\baselineskip = 10 truept
\centerline{\bf \linkint. Introduction}  
\bigskip

Since any three dimensional Einstein manifold has constant curvature, the
essential study of \Se manifolds begins in dimension five. Moreover, since
a complete \Se manifold is necessarily spin with finite fundamental group,
Smale's classification [Sm] of simply connected compact 5-manifolds with
spin applies. If there is no torsion in $H_2$ then Smale's theorem says that
any such 5-manifold is diffeomorphic to $S^5\#k(S^2\times S^3)$ for some
$k.$ In the classification of \Se manifolds it is judicious to distinguish
between regular \Se manifolds and non-regular ones. The compact simply
connected five dimensional manifolds admitting a regular \Se structure have
been classified by Friedrich and Kath [FK], and it follows from the
classification of smooth del Pezzo surfaces admitting K\"ahler-Einstein
metrics due to Tian and Yau [TY]. These 5-manifolds are $S^5$ and
$\#k(S^2\times S^3)$ for $k=1,3,\cdots,8.$ For $k=3,\cdots,8$ they are
are circle bundles over $\bbp^2$ blown up at $k$ generic points, whereas for
$k=1$ this is the homogeneous Stiefel manifold $V_{4,2}(\bbr)$ which is a
circle bundle over $\bbp^1\times \bbp^1.$ Notice that $k=2$ is missing from
this list as are the circle bundles over $\bbp^2$ blown-up at one or two
points. The reason for this is Matsushima's well-known obstruction to the
existence of K\"ahler-Einstein metrics when the complex automorphism group is
not reductive. Thus, any \Se structure on the connected sum of two copies of
$S^2\times S^3$ must be non-regular. In this note we prove the existence of
such a non-regular \Se metric thus filling this ``$k=2$ gap". 
It is an interesting
question as to whether \Se structures exist on $\#k(S^2\times S^3)$ for
$k>8.$  Of course, if such structures exist they must be non-regular. 

\noindent{\sc Note Added}: The ideas of this paper have been developed much 
further in [BGN1,BGN2] where the existence of infinite families of 
Sasakian-Einstein metrics on the k-fold connected sum $k\#(S^2\times S^3)$ is 
proven for $k=2,\cdots, 9.$ The moduli of such structures is also discussed. However, it
is still an open question as to whether there are \Se metrics on $\#k(S^2\times S^3)$  for 
$k>9.$

We also prove the existence of two inhomogeneous non-regular \Se metrics on
$S^2\times S^3.$ Recently it has been shown that the manifold $S^2\times S^3$
admits quite a few Einstein metrics. First Wang and Ziller [WZ] proved the
existence of a countable number of homogeneous Einstein metrics on $S^2\times
S^3.$ More recently B\"ohm [B\"oh] showed the existence of a countable
number of cohomogeneity one Einstein metrics on $S^2\times S^3.$ Our new
\Se metrics are also inhomogeneous and they are not isometric to any
of the B\"ohm's examples. Explicitly, we prove

\noindent{\sc Theorem A}: \tensl There exists a non-regular \Se metric on 
$(S^2\times S^3)\# (S^2\times S^3).$ \tenrm

\noindent{\sc Theorem B}: \tensl There exist
two inequivalent inhomogeneous \Se
metrics on $S^2\times S^3.$ These metrics are inequivalent as
Riemannian metrics to the inhomogeneous metrics of B\"ohm. Hence,
$S^2\times S^3$ admits at least three distinct \Se metrics.\tenrm

A Sasakian structure on a manifold defines several interesting objects. It
defines a one-dimensional foliation, a CR-structure, and a contact structure.
Indeed, combining all three of these, it defines a Pfaffian structure with a
transverse K\"ahler geometry. Now in each case above the \Se structure is
unique within the CR-structure. Thus, on $S^2\times S^3$ we have three
distinct CR-structures. It is interesting to ask the question as to whether
the three \Se structures on $S^2\times S^3$ belong to distinct contact
structures. This is a more subtle question as contact geometry has no local
invariants. Perhaps there is a connection between certain link invariants and
contact invariants as suggested by Arnold [Arn].

Our method of proofs of Theorems A and B is to consider the Sasakian geometry
of links of isolated hypersurface singularities defined by weighted homogeneous
polynomials. We then make use of a recent result of Demailly and Koll\'ar [DK]
proving the existence of K\"ahler-Einstein metrics on certain del Pezzo
orbifolds given as hypersurfaces in certain weighted projective spaces. The
links which can then be represented as the total space of V-bundles over these
orbifolds admit \Se metrics. We then use a well-known algorithm of Milnor and
Orlik [MO] to compute the characteristic polynomials of the monodromy maps
associated to the links. This allows us to determine the second Betti number
of the link. Then using a method of Randell [Ran] we can show that the links
have no torsion, and apply Smale's classification theorem. 

\noindent{\sc Acknowledgments}: We would like to thank Alex Buium and Michael
Nakamaye for helpful discussions. We also want to thank J\'anos Koll\'ar for
several valuable e-mail communications as well as his interest in our work.

\bigskip
\baselineskip = 10 truept
\centerline{\bf \linksas. The Sasakian Geometry of Links of Weighted
Homogeneous Polynomials}   \bigskip

In this section we discuss the Sasakian geometry of links of isolated
hypersurface singularities defined by weighted homogeneous polynomials.
Consider the affine space $\bbc^{n+1}$ together with a weighted
$\bbc^*$-action given by $(z_0,\ldots,z_n)\mapsto
(\grl^{w_0}z_0,\ldots,\grl^{w_n}z_n),$ where the {\it weights} $w_j$ are
positive integers. It is convenient to view the weights as the components of a
vector $\bfw\in (\bbz^+)^{n+1},$ and we shall assume that
$\gcd(w_0,\ldots,w_n)=1.$ Let $f$ be a quasi-homogeneous polynomial, that is
$f\in \bbc[z_0,\ldots,z_n]$ and satisfies
$$f(\grl^{w_0}z_0,\ldots,\grl^{w_n}z_n)=\grl^df(z_0,\ldots,z_n),
\leqno{\linksas.1}$$
where $d\in \bbz^+$ is the degree of $f.$ We are interested in the {\it
weighted affine cone} $C_f$ defined by
the equation $f(z_0,\ldots,z_n)=0.$ We shall assume that the origin in
$\bbc^{n+1}$ is an isolated singularity, in fact the only singularity, of
$f.$ Then the link $L_f$ defined by 
$$L_f= C_f\cap S^{2n+1}, \leqno{\linksas.2}$$
where 
$$S^{2n+1}=\{(z_0,\ldots,z_n)\in \bbc^{n+1}|\sum_{j=0}^n|z_j|^2=1\}$$
is the unit sphere in $\bbc^{n+1},$ is a smooth manifold of dimension $2n-1.$ 
Furthermore, it is well-known [Mil] that the link $L_f$ is $(n-2)$-connected.

On $S^{2n+1}$ there is a well-known [YK] ``weighted'' Sasakian structure  
$(\xi_\bfw,\eta_\bfw,\Phi_\bfw,g_\bfw)$ which in the standard coordinates
$\{z_j=x_j+iy_j\}_{j=0}^n$ on $\bbc^{n+1}=\bbr^{2n+2}$ is determined by
$$\eta_\bfw = {\sum_{i=0}^n(x_idy_i-y_idx_i)\over\sum_{i=0}^n
w_i(x_i^2+y_i^2)}, \qquad \xi_\bfw
=\sum_{i=0}^nw_i(x_i\partial_{y_i}-y_i\partial_{x_i}),$$
and the standard Sasakian structure $(\xi,\eta,\Phi,g)$ on $S^{2n+1}.$
Explicitly, we have
$$\eqalign{\Phi_\bfw &=\Phi-\Phi\xi_w\otimes \eta_w \cr g_\bfw &={1\over
\eta(\xi_\bfw)} [g-\eta_\bfw\otimes \xi_\bfw\rfloor g-
\xi_\bfw\rfloor g\otimes \eta_\bfw +g(\xi_\bfw,\xi_\bfw)\eta_\bfw\otimes
\eta_\bfw] +  \eta_\bfw\otimes \eta_\bfw.} \leqno{\linksas.3}$$

Now, by equation \linksas.1, the $\bbc^*(\bfw)$ action on $\bbc^{n+1}$ 
restricts
to an action on $C_f,$ and the associated $S^1$ action restricts to an action
on both $S^{2n+1}$ and $L_f.$ It follows that $\xi_\bfw$ is tangent to the
submanifold $L_f$ and, by abuse of notation, we shall denote by
$\xi_\bfw,\eta_\bfw,\Phi_\bfw,g_\bfw$ the corresponding tensor fields 
on both $S^{2n+1}$ and $L_f.$  Now  $\Phi_\bfw$ coincides with $\Phi$ on the
contact subbundle $\cald$ on $S^{2n+1}$ which defines an integrable almost
complex structure on $\cald.$ Moreover, since $f$ is a holomorphic function on
$\bbc^{n+1}$ the Cauchy-Riemann equations imply that for any smooth section
$X$ of $\cald$ we have $\Phi_\bfw X(f)=0.$ Thus, $L_f$ is an invariant
submanifold of $S^{2n+1}$ with its weighted Sasakian structure. We have
arrived at a theorem  given by Takahashi [Tak, YK] in the case of
Brieskorn-Pham links and we have seen that Takahashi's proof easily
generalizes to the case of arbitrary weighted homogeneous hypersurface
singularities. 

\noindent{\sc Theorem} \linksas.4: \tensl The quadruple
$(\xi_\bfw,\eta_\bfw,\Phi_\bfw,g_\bfw)$ gives $L_f$ a quasi-regular Sasakian
structure.  \tenrm

Actually as with K\"ahler structures there are many Sasakian structures on a
given Sasakian manifold. In fact there are many Sasakian structures 
which have $\xi$ as its characteristic vector field. To see this let
$(\xi,\eta,\Phi,g)$ be a Sasakian structure on a smooth manifold (orbifold)
$M,$ and consider a deformation of this structure by adding to $\eta$ a 
continuous one parameter family of 1-forms $\grz_t$ that are basic with
respect to the characteristic foliation. We require that the 1-form
$\eta_t=\eta +\grz_t$  satisfy the conditions  
$$\eta_0=\eta, \qquad \grz_0=0,\qquad \eta_t\wedge (d\eta_t)^n\neq 0~~
\forall~~ t\in [0,1]. \leqno{\linksas.5}$$  
This last non-degeneracy condition implies that $\eta_t$ is a contact form on
$M$ for all $t\in [0,1].$ Then by Gray's Stability Theorem [MS] $\eta_t$
belongs to the same contact structure as $\eta.$ Moreover, since $\grz_t$ is
basic $\xi$ is the Reeb (characteristic) vector field associated to $\eta_t$
for all $t.$ Now let us define
$$\eqalign{\Phi_t&=\Phi -\xi\otimes \grz_t\circ \Phi \cr
           g_t&=g+d\grz_t\circ (\Phi\otimes \hbox{id}) + \grz_t\otimes\eta
+\eta\otimes \grz_t +\grz_t\otimes \grz_t.} \leqno{\linksas.6}$$
Note that it is not at all clear from this definition that $g_t$ is a
Riemannian metric, but we shall check this below. We have

\noindent{\sc Theorem} \linksas.7: \tensl Let $(M,\xi,\eta,\Phi,g)$ be a
Sasakian manifold. Then for all $t\in [0,1]$ and every basic 1-form $\grz_t$
such that $d\grz_t$ is of type $(1,1)$  and such that \linksas.5 holds
$(\xi,\eta_t,\Phi_t,g_t)$ defines a Sasakian structure on $M$ belonging to the
same underlying contact structure as $\eta.$ \tenrm  

\noindent{\sc Proof}: The conditions of \linksas.5 guarantee that
$(\xi,\eta_t,\Phi_t,g_t)$ defines a Pfaffian structure, i.e. a contact
structure with a fixed contact 1-form. We need to check that it is a metric
contact structure and that it is normal. It is easy to check that
the metric $g_t$ of \linksas.6 can be rewritten as
$$g_t=d\eta_t\circ(\Phi_t\otimes \hbox{id})+\eta_t\otimes \eta_t.
\leqno{\linksas.8}$$
It follows from the fact $d\eta_t$ is type $(1,1)$ on the contact bundle
$\cald_t=\hbox{ker}~\eta_t$ that $g_t$ is a symmetric bilinear form and then a
straightforward computation checks the compatibility condition 
$$g_t(\Phi_tX,\Phi_tY)=g_t(X,Y)-\eta_t(X)\eta_t(Y).$$
The positive definiteness of $g_t$ follows from the positive definiteness of
$g$ and the non-degeneracy condition in \linksas.5. Moreover, one easily checks
the identity $\Phi_t^2=-\hbox{id}+\xi\otimes \eta_t.$ Next we check normality
which amounts to checking two conditions, that the almost CR structure defined
by $\Phi_t$ on the contact bundle $\cald_t$ is integrable, and that $\xi$ is a
Killing vector field for the metric $g_t.$ The last condition is equivalent to
vanishing of the Lie derivative $\call_\xi\Phi_t$ for all $t$ which follows
immediately from the first of equations \linksas.6 and the facts that it holds
for $t=0$ and that $\grz_t$ is basic. Integrability follows from the fact that
the almost CR structure defined by $\Phi$ on $\cald$ is integrable, and that
the first of equations \linksas.6 is just the projection of the image of
$\Phi$ onto $\cald_t.$ \hfill\za

In general these structures are inequivalent and the moduli
space of Sasakian structures having the same characteristic vector field
is infinite dimensional. Indeed since the link of a hypersurface is determined
by the $S^1(\bfw)$ action we formulate the following:

\noindent{\sc Definition} \linksas.8: \tensl A Sasakian structure
$(\xi,\eta,\Phi,g)$ on $L_f$ is said to be {\it compatible} with the link
$L_f$ if $\xi$ is a generator of the $S^1$ action on $L_f.$  We say
that $\xi$ is the {\it standard} generator if $\xi=\xi_\bfw$, where $\bfw$ is
the weight vector of $L_f$ satisfying $\gcd(w_0,\ldots,w_n)=1.$  \tenrm  

For every compatible Sasakian structure there is one with a standard
generator, and hereafter we shall always choose the standard generator for 
a compatible Sasakian structure unless otherwise stated. We are interested in
the following question:

\noindent{\sc Problem} \linksas.9: \tensl Given a link $L_f$ with a given
Sasakian structure $(\xi,\eta,\Phi,g),$ when can we find a 1-form $\grz$ such
that the deformed structure $(\xi,\eta+\grz,\Phi',g')$ is \Se? \tenrm

This is a Sasakian version of the Calabi problem for the link $L_f$ which is 
discussed further in [BGN1].  Here we use the fact that the leaf space of the
characteristic foliation of a Sasakian structure of a link $L_f$ is a compact
K\"ahler orbifold together with recent results of Demailly and Koll\'ar [DK] on
the existence of K\"ahler-Einstein orbifold metrics on certain singular del
Pezzo surfaces to construct the \Se metrics on the corresponding link.

\bigskip
\baselineskip = 10 truept
\centerline{\bf \linkke. K\"ahler-Einstein Orbifolds and \Se Manifolds}
\bigskip

Given a sequence $\bfw =(w_0,\ldots,w_n)$ of positive integers one can form
the graded polynomial ring $S(\bfw)=\bbc[z_0,\ldots,z_n]$, where $z_i$ has
grading or {\it weight} $w_i.$ The weighted projective space [BR, Del, Dol, 
Fle]
$\bbp(\bfw)=\bbp(w_0,\ldots,w_n)$ is defined to be the scheme
$\hbox{Proj}(S(\bfw)).$  It is the quotient space 
$(\bbc^{n+1}-\{0\})/\bbc^*(\bfw)$, where $\bbc^*(\bfw)$ is the weighted action
defined in section \linksas. Clearly, $\bbp(\bfw)$ is also the quotient of
the weighted Sasakian sphere
$S_\bfw^{2n+1}=(S^{2n+1},\xi_\bfw,\eta_\bfw,\Phi_\bfw,g_\bfw)$ by the
weighted circle action $S^1(\bfw)$ generated by $\xi_\bfw.$ As such
$\bbp(\bfw)$ is also a compact complex orbifold with an induced K\"ahler
structure. At times it will be important to distinguish between $\bbp(\bfw)$
as a complex orbifold and $\bbp(\bfw)$ as an algebraic variety.  

Now the cone $C_f$ in $\bbc^{n+1}$ cuts out a hypersurface $\calz_f$ of
$\bbp(\bfw)$ which is also a compact orbifold with an induced K\"ahler
structure  $\gro_\bfw.$ So there is a commutative diagram
$$\matrix{L_f &\ra{2.5}& S^{2n+1}_\bfw&\cr
  \decdnar{\pi}&&\decdnar{} &\cr
   \calz_f &\ra{2.5} &\bbp(\bfw),&\cr}\leqno{\linkke.1}$$
where the horizontal arrows are Sasakian and K\"ahlerian embeddings,
respectively, and the vertical arrows are orbifold Riemannian submersions.
Furthermore, by the inversion theorem of [BG1] $L_f$ is the total space of the
principal $S^1$ V-bundle over the orbifold $\calz_f$ whose first Chern class is
$[\gro_\bfw]\in H^2_{orb}(\calz_f,\bbz),$ and $\eta_\bfw$ is the connection in
this V-bundle whose curvature is $\pi^*\gro_\bfw.$ For further discussion of
the orbifold cohomology groups we refer the reader to [Hae] and [BG1].

Now $H^2_{orb}(\calz_f,\bbz)\otimes \bbq\approx H^2(\calz_f,\bbq),$ so
$[\gro_\bfw]$ defines a rational class in the usual cohomology.
Thus, by Baily's [Bai] projective embedding theorem for orbifolds, $\calz_f$
is a projective algebraic variety with at most quotient singularities. It
follows that as an algebraic variety $\calz_f$ is normal. As with $\bbp(\bfw)$
it is important to distinguish between $\calz_f$ as an orbifold and
$\calz_f$ as an algebraic variety. For example the singular loci may differ.
There are examples where the orbifold singular locus $\grS^{orb}(\calz_f)$ has
codimension one over $\bbc,$ whereas $\grS^{alg}(\calz_f)$ always has
codimension $\geq 2$ by normality. Indeed there are examples with
$\grS^{orb}(\calz_f)\neq \emptyset,$ but $\calz_f$ is smooth as an algebraic
variety. In these cases an orbifold metric is NOT a Riemannian metric in the
usual sense.

We are interested in finding \Se structures $(\xi_\bfw,\eta'_\bfw,\Phi'_\bfw,
g'_\bfw)$ on $L_f$ with the same Sasakian vector field $\xi_\bfw$ as the
original Sasakian structure, that is \Se structures that are compatible with
the link $L_f.$ Since a \Se metric necessarily has positive Ricci curvature,
we see that a necessary condition for a \Se metric $g$ to exist on $L_f$ is
that there be an orbifold metric $h$ on $\calz_f$ with positive Ricci form
$\grr_h$ such that $g=\pi^*h+\eta\otimes \eta.$  The positive definiteness of
$\grr_h$ follows from the relation
$$\pi^*\Ric_h=\Ric_g|_{\cald\times \cald} + 2g|_{\cald\times \cald}. 
\leqno{\linkke.2}$$
But since $\grr_h$ represents the first Chern class $c_1(K^{-1})$ of the
anti-canonical line  V-bundle $K^{-1}$ we see that $\calz_f$ must be a Fano
orbifold, i.e., some power of $K^{-1}$ is invertible and ample. Here $K$ is the
canonical line V-bundle of Baily [Bai].  By the previous section and the
well-known theory of \Se structures [BG1] we have

\noindent{\sc Theorem} \linkke.3: \tensl The link $L_f$ has a compatible
\Se structure if and only if the Fano orbifold $\calz_f$ admits a compatible
K\"ahler-Einstein orbifold metric of scalar curvature $4n(n+1).$ \tenrm

Of course, if we find a K\"ahler-Einstein metric of positive scalar curvature
on the orbifold $\calz_f,$ we can always rescale the metric so that the scalar
curvature in $4n(n+1).$

We want conditions that guarantee that the hypersurfaces be Fano, but first we
restrict somewhat the hypersurfaces that we treat.   

\noindent{\sc Definition} \linkke.3 [Fle]: \tensl 
\item{(1)} A  weighted projective
space $\bbp(w_0,\ldots,w_n)$ is said to be {\it well-formed} if 
$$\gcd(w_0,\ldots,\hat{w_i},\ldots,w_n)=1~\hbox{for all $i=1,\ldots,n.$}$$

\item{(2)} A hypersurface in a weighted projective space $\bbp(w_0,\ldots,w_n)$
is said to be {\it well-formed} if in addition it contains no codimension 2
singular stratum of $\bbp(w_0,\ldots,w_n),$

\noindent where the hat means delete that element. \tenrm

In (1) of definition \linkke.3, we prefer the terminology introduced in [BR]
for an equivalent condition which seems to have first appeared in [Dol]. So if
a weighted projective space satisfies (1) of definition \linkke.3 we say that
$\bfw$ is {\it normalized}. In [Fle] it is shown that (2) of definition
\linkke.3 can be formulated solely in terms of the weighted homogeneous
polynomial $f.$ That is, a general hypersurface in $\bbp(\bfw)$ defined by a
weighted homogeneous polynomial $f$ is well-formed if and only if
$\bfw$ is normalized and 
$$\gcd(w_0,\ldots,\hat{w_i},\ldots,\hat{w_j},\ldots,w_n)~|~d \qquad \hbox{for
all distinct $i,j=0,\ldots, n.$} \leqno{\linkke.4}$$ 
In this case we shall also say that the weighted homogeneous polynomial $f$
defining the hypersurface is well-formed. The following proposition is
essentially an exercise in [Kol]:

\noindent{\sc Proposition} \linkke.5: \tensl Let $f_\bfw$ be a well-formed
weighted homogeneous polynomial with weights $\bfw$ and degree $d.$ Then the
hypersurface $\calz_f$ in $\bbp(w_0,\ldots,w_n)$ is Fano if and only if the
following condition holds
\item{$(\ast)$} $d< |\bfw|=\sum_{i=0}^nw_i.$

\noindent Thus if $\bfw$ is normalized and satisfies \linkke.4, a necessary
condition for the link $L_f$ to admit a compatible \Se metric is that $(\ast)$
be satisfied.  \tenrm

\noindent{\sc Proof}: Following [Mo,BR] in any weighted projective space
$\bbp(\bfw)$ we define the Mori-singular locus as follows: Let
$m=\lcm(w_0,\ldots,w_n)$ and consider for each prime divisor $p$ of $m$ the
subscheme $S_p$ cut out by the ideal $I_p$ generated by those indeterminates
$z_i$ such that $p\not| w_i,$ and define the Mori-singular locus by
$$S(\bfw)=\cup_{p|m} S_p,$$ 
and the Mori-regular locus by
$$\bbp^0(\bfw)=\bbp(\bfw)-S(\bfw).$$ 
Now $\bbp^0(\bfw)$ is contained in the
smooth locus of $\bbp(\bfw)$ but can be a proper subset of it. Thus, $S(\bfw)$
contains the singular locus $\grS(\bfw)$ of $\bbp(\bfw)$ but in general is
larger. However, when $\bfw$ is normalized the local uniformizing groups of the
orbifold $\bbp(\bfw)$ contain no quasi-reflections as in [Fle]. It follows in
this case that the singular locus $\grS(\bfw)$ coincides with the
Mori-singular locus $S(\bfw),$ and $\bbp^0(\bfw)$ is precisely the smooth
locus of $\bbp(\bfw).$ 

Now the sheaves $\calo_{\bbp(\bfw)}(n)$ for $n\in \bbz$ are not in general
invertible, but they are reflexive [BR]. Moreover, the Mori-regular locus
$\bbp^0(\bfw)$ is the largest open subset $U$ of $\bbp(\bfw)$ on which
\item{(1)} $\calo_{\bbp(\bfw)}(1)|U$ is invertible, and 
\item{(2)} the natural map $\calo_{\bbp(\bfw)}(1)^{\otimes n}|U\ra{1.4}
\calo_{\bbp(\bfw)}(n)|U$ is an isomorphism for each $n\in \bbz^+$ [Mo].

\noindent Furthermore, if $\bfw$ is normalized $\bbp^0(\bfw)$ is characterized
by condition (1) alone.  As an orbifold $\bbp(\bfw)$ has a canonical V-bundle
$K_{\bbp(\bfw)}$ whose associated sheaf, the dualizing sheaf
$\gro_{\bbp(\bfw)}$, is isomorphic to $\calo_{\bbp(\bfw)}(-|\bfw|).$ Generally
the adjunction formula does not hold, but if $\bfw$ is normalized we have
$$\gro_{\calz_f}\cong \calo_{\calz_f}(d-|\bfw|). \leqno{\linkke.6}$$
These sheaves are invertible on $\bbp^0(\bfw)\cap \calz_f$ and we have 
$$\gro^{-1}_{\calz_f\cap \bbp^0(\bfw)}\cong \calo_{\calz_f\cap
\bbp^0(\bfw)}(|\bfw|-d).$$
Moreover, since $\calz_f$ is well-formed, $\calz_f\cap S(\bfw)=\calz_f\cap
\grS(\bfw)$ has codimension $\geq 2$ in $\calz_f.$ It follows by a standard
result (cf. [KMM] Lemma 0-1-10) that 
$$\gro^{-1}_{\calz_f}= \gri_*(\gro^{-1}_{\calz_f\cap \bbp^0(\bfw)}),\qquad
\calo_{\calz_f}(|\bfw|-d)= \gri_*(\calo_{\calz_f\cap
\bbp^0(\bfw)}(|\bfw|-d)),$$
where $\gri:\calz_f\cap \bbp^0(\bfw)\ra{1.3} \calz_f$ is the natural
inclusion. We then have   
$$\gro^{-1}_{\calz_f}\cong \calo_{\calz_f}(|\bfw|-d).
\leqno{\linkke.7}$$ 
This proves the result. \hfill\za

There are further obstructions to the existence of a \Se structure on $L_f,$
namely those coming from the well-known failure of the zeroth order a priori
estimate for solving the Monge-Ampere equation for the Calabi problem on
$\calz_f.$ Now the link $L_f$ admits a compatible \Se structure if and only if
$\calz_f$ has a K\"ahler-Einstein structure $\gro'$ in the same cohomology
class as $\gro_\bfw.$ So the well-known obstructions in the K\"ahler-Einstein
case, such as the non-reductiveness of the connected component of the complex
automorphism group, or the vanishing of the Futaki
invariant become obstructions to the existence of \Se metrics. The problem of
solving the Calabi problem in the positive case has drawn much attention in
the last decade [Sui, Ti1, Ti2, Ti3, TY, DT, Nad]
but still remains open. Here it
suffices to consider three examples given recently by Demailly and Kollar
[DK]. 

\bigskip
\baselineskip = 10 truept
\centerline{\bf \linkex. The Demailly-Koll\'ar Examples} 
\bigskip
 
In a recent preprint Demailly and Koll\'ar [DK] give a new derivation of
Nadel's existence criteria for positive K\"ahler-Einstein metrics [Nad] which
is valid for orbifolds. Furthermore their method of implementing Nadel's
theorem does not depend on the existence of a large finite group of
symmetries, but rather uses intersection inequalities. In this sense it is
complementary to the work of Nadel [Nad] and others [Siu, Ti1-Ti4]. As an
application of their method Demailly and Koll\'ar [DK] construct three new del
Pezzo orbifolds which admit K\"ahler-Einstein metrics. Explicitly, they prove

\noindent{\sc Theorem} \linkex.1 [DK]: \tensl Let $\calz_f$ be a del Pezzo
orbifold, and let $f$ be given by one of the following three quasi-homogeneous
polynomials of degree $d$ whose zero set $\calz_f$ is a surface in the weighted
projective space $\bbp(w_0,w_1,w_2,w_3):$
\item{(1)} $f=z_0^5z_1+z_0z_2^3+z_1^4+z_3^3$ with $\bfw=(9,15,17,20)$ and
$d=60.$ 
 \item{(2)} $f=z_0^{17}z_2+z_0z_1^5+z_1z_2^3+z_3^2$ with $\bfw=(11,49,69,128)$ 
and $d=256.$ 
\item{(3)} $f=z_0^{17}z_1+z_0z_2^3+z_1^5z_2+z_3^2$ with $\bfw=(13,35,81,128)$
and $d=256.$

\noindent Then $\calz_f$ admits a K\"ahler-Einstein orbifold metric.\tenrm

We denote these del Pezzo surfaces by
$\calz_{60},\calz_{256}^{(1)},\calz_{256}^{(2)},$ respectively. Actually more
is true. By a result of Bando and Mabuchi [BM] (which also holds in the case of
orbifolds) the K\"ahler-Einstein metric is unique up to complex automorphisms.
Let us denote the corresponding links $L_f$ of these del Pezzo surfaces by
$L_{60},L_{256}^{(1)},L_{256}^{(1)},$ respectively. Due to the form of the
polynomials as $f=g+z_3^m,$ there is a well-known description of both $\calz_f$
and $L_f$ in terms of branched covers (cf. [DuKa]). We briefly describe this
geometry here. $\calz_{60}$ is a 3-fold cover of $\bbp(9,15,17)$ branched
over the curve $C_{60}=\{z_0^5z_1+z_0z_2^3+z_1^4=0\}.$ Similarly, the surfaces
$\calz_{256}^{(i)}$ are 2-fold covers of $\bbp(11,49,69)$ and
$\bbp(13,35,81),$ respectively, branched over the curves 
$C_{256}^{(1)}=\{z_0^{17}z_2+z_0z_1^5+z_1z_2^3=0\},$ and
$C_{256}^{(2)}=\{z_0^{17}z_1+z_0z_2^3+z_1^5z_2=0\},$ respectively. In 
section \linktop~ we show that in each case the genus of these curves is zero,
so they are all $\bbp^1$'s (there are no quotient singularities in dimension
one). Similarly, $L_{60}$ is a 3-fold cyclic cover of $S^5$ branched over the
Seifert manifold $K(4,3,5;IV),$ and $L_{256}^{(i)}$ are double covers of $S^5$
branched over the Seifert manifolds $K(17,3,5;V)$ and $K(17,5,3;V)$ for
$i=1,2,$ respectively. The notation $K(a_0,a_1,a_2;\cdot)$ is that of Orlik
[Or1]. These Seifert manifolds are quite complicated. For example, they all
have infinite fundamental group, but have finite abelianization $H_1,$ and
they all are Seifert fibrations over the Riemann sphere. As we shall see
shortly the links $L_{256}^{(1)}$ and $L_{256}^{(2)}$ have the same
characteristic polynomial although their Sasakian structures are distinct
since their leaf holonomy groups are different. In particular they can be
distinguished by their orders [BG1]. One easily finds by analyzing the
orbifold singularity structure that $\hbox{ord}(L_{256}^{(1)})=37191=3\cdot
7^2\cdot11\cdot 23,$ while
$\hbox{ord}(L_{256}^{(2)})=36855=3^4\cdot5\cdot7\cdot 13.$ Now combining the
Demailly-Koll\'ar Theorem with Theorem \linkke.3 gives 

\noindent{\sc Theorem} \linkex.2: \tensl The simply connected 5-manifolds
$L_{60},L_{256}^{(i)}$ admit compatible \Se metrics. Furthermore, these \Se
metrics are unique within the CR structure up to a CR automorphism. \tenrm

\bigskip
\baselineskip = 10 truept
\centerline{\bf \linkmi. Link Invariants and the Milnor Fibration} 
\bigskip

Recall the well-known construction of Milnor [Mil] for isolated hypersurface
singularities: There is a fibration of $(S^{2n+1}-L_f)\ra{1.3} S^1$ whose fiber
$F$ is an open manifold that is homotopy equivalent to a bouquet of n-spheres
$S^n\vee S^n\cdots \vee S^n.$ The {\it Milnor number} $\mu$ of $L_f$ is the
number of $S^n$'s in the bouquet. It is an invariant of the link which
can be calculated explicitly in terms of the degree $d$ and weights
$(w_0,\ldots,w_n)$ by the formula [MO]
$$\mu =\mu(L_f)=\prod_{i=0}^n\bigl({d\over w_i}-1\bigr).\leqno{\linkmi.1}$$
One immediately has

\noindent{\sc Proposition} \linkmi.2: \tensl The Milnor numbers of the 
simply connected \Se 5-manifolds $L_{60},L_{256}^{(1)},L_{256}^{(2)}$ are
$$\mu(L_{60})=86,\qquad \mu(L_{256}^{(i)})= 255.$$ \tenrm

The closure $\bar{F}$ of $F$ has the same homotopy type as $F$ and is a compact
manifold with boundary precisely the link $L_f.$ So the reduced homology of
$F$ and $\bar{F}$ is only non-zero in dimension $n$ and $H_n(F,\bbz)\approx
\bbz^\mu.$ Using the Wang sequence of the Milnor fibration together with
Alexander-Poincare duality gives the exact sequence [Mil]
$$0\ra{1.5} H_n(L_f,\bbz)\ra{1.5} H_n(F,\bbz)
\fract{\bbi -h_*}{\ra{1.5}} H_n(F,\bbz) \ra{1.5} H_{n-1}(L_f,\bbz)\ra{1.5} 0,
\leqno{\linkmi.3}$$ 
where $h_*$ is the {\it monodromy} map (or characteristic
map) induced by the $S^1_\bfw$ action. From this we see that
$H_n(L_f,\bbz)=\ker(\bbi -h_*)$ is a free Abelian group, and $H_{n-1}(L_f,\bbz)
=\coker(\bbi -h_*)$ which in general has torsion, but whose free part equals 
$\ker(\bbi -h_*).$ So the topology of $L_f$ is encoded in the monodromy map
$h_*.$  There is a well-known algorithm due to Milnor and Orlik [MO]
for computing the free part of $H_{n-1}(L_f,\bbz)$ in terms of the
characteristic polynomial $\grD(t)=\det(t\bbi -h_*),$ namely the Betti number
$b_n(L_f)=b_{n-1}(L_f)$ equals the number of factors of $(t-1)$ in $\grD(t).$ 
We now compute the characteristic polynomials $\grD(t)$ for our examples.

\noindent{Proposition} \linkmi.4: \tensl The characteristic polynomials
$\grD(t)$ of the simply connected \Se 5-manifolds
$L_{60},L_{256}^{(1)},L_{256}^{(2)}$ are given by 
$$\matrix{\grD(t)= &(t-1)^2(t^{14}+t^{13}+\cdots
+1)(t^{15}+1)(t^{30}+1)(t^4+t^3+t^2+t+1)\cr &\times
(t^{5}+1)(t^{10}+1)(t^4-t^2+1)(t^2-t+1)}$$
for $L_{60}$ and
$$\grD(t)=(t-1)(t^2+1)(t^4+1)(t^8+1)(t^{16}+1)(t^{32}+1)(t^{64}+1)(t^{128}+1)$$
for $L_{256}^{(1)}$ and $L_{256}^{(2)}.$ \tenrm

\noindent{\sc Proof}: The Milnor and Orlik [MO] algorithm for computing the
characteristic polynomial of the monodromy operator for weighted homogeneous
polynomials is as follows: First associate to any monic polynomial $F$ with
roots $\gra_1,\ldots,\gra_k\in \bbc^*$ its divisor 
$$\div F= <\gra_1>+\cdots+<\gra_k>$$
as an element of the integral ring $\bbz[\bbc^*]$ and let $\grL_n= \div
(t^n-1).$  The rational weights $w'_i$ used in [MO] are related to our integer
weights $w_i$ by $w_i'={d\over w_i},$ and we write the $w'_i={u_i\over v_i}$ in
irreducible form. So for the surface of degree $60$ we have 
$$w'_0={20\over 3},\quad w'_1=4, \quad w'_2={60\over 17}, \quad  w'_3=3.$$
Then by Theorem 4 of [MO] the divisor of the characteristic polynomial is 
$$\div \grD(t)=({\grL_{20}\over 3}-1)(\grL_4-1)({\grL_{60}\over
17}-1)(\grL_3-1).$$
Using the relations $\grL_a\grL_b=\gcd(a,b)\grL_{lcm(a,b)}$ we find in this
case
$$\div 
\grD(t)=\grL_{60}+\grL_{20}+\grL_{12}-\grL_4-\grL_3+1.\leqno{\linkmi.5}$$  
Then $\grD(t)$ is obtained from the formula  $$\div
\grD(t)= 1-\sum {s_j\over j}\grL_j$$ yielding
$$\grD(t)={(t-1)\over \prod(t^j-1)^{s_j/j}}.$$
We see from \linkmi.5 that the only nonzero $s_j$'s are 
$$s_{60}=-60,\quad s_{20}=-20, \quad s_{12}=-12,\quad s_4=4, \quad s_3=3,$$
and $\grD(t)$ becomes
$$\grD(t)={(t-1)(t^{60}-1)(t^{20}-1)(t^{12}-1)\over (t^4-1)(t^3-1)}$$
which we see easily reduces the the expression given.

Similar computations for the two surfaces of degree $256$ show that 
$$\div \grD(t)= \grL_{256}-\grL_2+1, \leqno{\linkmi.6}$$
giving the characteristic polynomial noted above. \hfill\za

\bigskip
\baselineskip = 10 truept
\centerline{\bf \linktop. The Topology of $L_f$} 
\bigskip

Since the multiplicity of the root $1$ in $\grD(t)$ is precisely the second
Betti number of $L_f$ we have an immediate corollary of Proposition \linkmi.4: 

\noindent{\sc Corollary} \linktop.1: \tensl  The second Betti numbers of the
5-manifolds $L_{60},L_{256}^{(i)}$ are $b_2(L_{60})=2,$ and
$b_2(L_{256}^{(i)})=1$ for $i=1,2.$ \tenrm

One can give another proof of this corollary by making use of a formula due to
Steenbrink [Dim]. One considers the Milnor algebra
$$M(f)= {\bbc[z_0,z_1,z_2,z_3]\over \theta(f)}, \leqno{\linktop.2}$$
where $\theta(f)$ is the Jacobi ideal. $M(f)$ has a natural $\bbz$ grading 
$$M(f)=\oplus_{i\geq 0}M(f)_i,$$ and Steenbrink's formula determines the
primitive Hodge numbers  
$$h^{i,n-i-1}_0= \dim_\bbc M(f)_{(i+1)d-w} \leqno{\linktop.3}$$
of the projective surface $f=0$ in the weighted projective space $\bbp(\bfw).$
Then $b_2(L_f)=h^{1,1}_0$ can be computed by finding a basis of residue
classes in $M(f)_{2d-|\bfw|}.$ Note also that since $\calz_f$ is a del Pezzo surface,
we have $h^{2,0}=h^{0,2}=0.$

Steenbrink's formula can also be used to compute the Hirzebruch signature
$\grt$ of the orbifolds $\calz_f$ as well as the genus of the curves $C_{60}$
and $C_{256}^{(i)}$ discussed at the end of the last section. We have

\noindent{\sc Proposition} \linktop.4: \tensl The Fano orbifold $\calz_{60}$
has $\grt=-1$ while the orbifolds $\calz_{256}^{(i)}$ have $\grt=0.$ The curves
$C_{60}$ and $C_{256}^{(i)}$ are all isomorphic to the projective line
$\bbp^1.$  \tenrm

\noindent{\sc Proof}: The formula for the signature is [Dim]
$$\grt = 1+2\dim_\bbc M(f)_{d-|\bfw|} -\dim_\bbc M(f)_{2d-|\bfw|},$$
whereas the genus of the curves is given by  
$$g=h^{0,1}_0= \dim_\bbc M(f_0)_{d-|\bfw|},$$ 
where $f_0$ is a weighted homogeneous polynomial in $z_0,z_1,z_2$ which is
related to $f$ by $f=f_0+z_3^{a_3}.$ Now $d-|\bfw|=-1$ for all the surfaces and
$\dim_\bbc M(f)_{2d-|\bfw|}=h^{1,1}_0=b_2(L_f)$ which is $2$ for $\calz_{60}$
and $1$ for $\calz_{256}^{(i)}.$ Now for the curve $C_{60}$ we have
$d-|\bfw|=19,$ and one easily sees that $g= \dim_\bbc M(f_0)_{19}=0.$ For the
curves $C_{256}^{(i)}$ we have $d-|\bfw|=127,$ so $g= \dim_\bbc
M(f_0)_{127}.$ For example, for $C_{256}^{(1)}$ one easily checks that there
are no monomials of the form $z_0^az_1^bz_2^c$ such that $11a+49b+69c=127.$ In
both cases we find $g=0.$  \hfill\za

Next we turn to the calculation of the torsion in $H_2(L_f,\bbz).$ Actually we
show that it is zero. We follow a method due to Randell [Ran] for computing
the torsion of generalized Brieskorn manifolds (complete intersections of
Brieskorn's) where he verified a conjecture due to Orlik [Or2] for this class
of manifolds. Randell's methods apply to our more general weighted homogeneous
polynomials case. For any $x\in L_f$ let $\grG_x$ denote the isotropy subgroup
of the $S^1(\bfw)$-action on $L_f$ and by $|\grG_x|$ its order. Following
[Ran] for each prime $p$ we define the {\it $p$-singular set} by
$$S_p=\{x\in L_f|~p~\hbox{divides}~ |\grG_x|\}. \leqno{\linktop.5}$$  
Then we have

\noindent{\sc Lemma} \linktop.6: \tensl  Suppose that for each prime $p$ the
$p$-singular set $S_p$ is contained in a submanifold $S$ of codimension $4$ in
$L_f$ that is cut out by a hyperplane section of $\calz_f.$ Then the isomorphism holds: 
$$\hbox{Tor}(H_{n-1}(L_f,\bbz))\approx
\hbox{Tor}(H_{n-3}(S,\bbz)). \leqno{\linktop.7}$$ 
In particular, if the corresponding projective hypersurface $\calz_f$ is
well-formed, the isomorphism \linktop.7 holds. \tenrm

\noindent{\sc Proof}: The first statement is due to Randell [Ran]. To prove
the second statement we notice that the set $\cup_p S_p,$ where the union is
over all primes, is precisely the subset of $L_f$, where the leaf holonomy
groups are non-trivial and its projection to $\calz_f$ is just the orbifold
singular locus $\grS^{orb}(\calz_f).$ Now $\calz_f\subset \bbp(\bfw)$ and
$\bfw$ is normalized, so the fact that $\calz_f$ is well-formed says that the
orbifold singular locus $\grS^{orb}(\calz_f)$ has complex codimension at
least two in $\calz_f;$ hence, for each prime $p,$ $S_p$ has real codimension 
at least four in $L_f.$ \hfill\za

For the case at hand, $n=3$ so Randell's lemma says that $H_2(L_f,\bbz)$ is
torsion free. Thus, we have arrived at:

\noindent{\sc Lemma} \linktop.8: \tensl Let $L_f\subset \bbc^4$ be the link
of an isolated singularity defined by a weighted homogeneous polynomial $f$ in
four complex variables. Suppose further that the del Pezzo surface
$\calz_f\subset \bbp(\bfw)$ is well-formed, then
$\hbox{Tor}~(H_2(L_f,\bbz))=0.$ \tenrm    

Now a well-known theorem of Smale [Sm] says that any simply connected compact
5-manifold which is spin, and whose second homology group is torsion free, is
diffeomorphic to $S^5\#k (S^2\times S^3)$ for some non-negative integer $k.$
Furthermore, it is known [BG2,Mor] that any simply connected \Se manifold is
spin. Combining this with the development above gives

\noindent{\sc Theorem} \linktop.9: \tensl Let $L_f$ be the link associated to
a well-formed weighted homogeneous polynomial $f$ in four complex variables.
Suppose also that $L_f$ admits a \Se metric. Then $L_f$ is diffeomorphic to 
$S^5\#k (S^2\times S^3)$, where $k$ is the multiplicity of the
root $1$ of the characteristic polynomial $\grD(t)$ of $L_f.$ \tenrm

Let us now consider the links $L_{60},L_{256}^{(i)}$ of Theorem \linkex.2
which admit \Se metrics.

\noindent{\sc Theorem} \linktop.10: \tensl The link $L_{60}$ is diffeomorphic
to $(S^2\times S^3)\#(S^2\times S^3)$ while the links $L_{256}^{(i)}$ for
$i=1,2$ are diffeomorphic to the Stiefel manifold $S^2\times S^3.$ In
particular, $S^2\times S^3$ admits 3 distinct \Se structures. \tenrm   

\noindent{\sc Proof}: This follows from Smale's Theorem [Sm], Theorems
\linkex.2, \linktop.9, and Corollary \linktop.1 as soon as we check that the
weighted homogeneous polynomials are well-formed. But this follows easily from
the definition and equation \linkke.4. \hfill\za

\bigskip
\baselineskip = 10 truept
\centerline{\bf \linkpr. Proofs of the Main Theorems and Further Discussion} 
\bigskip

Theorem A of the Introduction now follows immediately from Theorems \linkex.2
and \linktop.10. Similarly for Theorem B these two theorems give the existence
of \Se metrics on $S^2\times S^3.$ Furthermore, since the Sasakian structures
are non-regular the metrics are inhomogeneous. The statement that the \Se
structures are inequivalent follows from the fact that their
characteristic foliations are inequivalent (indeed, they have different
orders). It remains to show that our \Se metrics are inequivalent as
Riemannian metrics to any of the metrics of B\"ohm. To see this we notice that
for all of B\"ohm's metrics on $S^2\times S^3,$ the connected component of the
isometry group $\gI_0(g)$ is $SO(3)\times SO(3)$ (see [B\"oh] or 
[Wa] Theorem 2.16). But
a Theorem of Tanno [Tan] says that the connected component of the group of
Sasakian automorphisms, $\gA_0(g)$ coincides with the connected component of
the group of isometries $\gI_0(g),$ and $\gA_0(g)$ has the form $\gG\times 
S^1$,
where the $S^1$ is generated by the Sasakian vector field $\xi.$ Thus, our
metrics are not isometric to the any of B\"ohm's metrics. \hfill\za

Smale's Theorem actually says more than we have mentioned. It says that the
torsion in $H_2$ must be of the form 
$$\hbox{Tor}(H_2(M,\bbz))\approx \bigoplus (\bbz_{q_i}\oplus \bbz_{q_i}).
\leqno{\linktop.11}$$
Thus we have

\noindent{\sc Theorem} \linkpr.1: \tensl Let $M$ be a complete simply
connected \Se 5-manifold. Then the torsion group $\hbox{Tor}(H_2(M,\bbz))$ 
must be of the form \linktop.11. In particular, the number of torsion
generators must be even. \tenrm

It is a very interesting question whether the form of the torsion actually
provides an obstruction to the existence of a \Se structure or whether it is
always satisfied. The general form of Randell's proof suggests that the
torsion may always vanish, but we do not yet have a proof of this even in the
case of hypersurfaces defined by weighted homogeneous polynomials when we drop
the well-formed assumption. If there were an honest obstruction, it would
provide a new type of obstruction to the existence of positive
K\"ahler-Einstein metrics on del Pezzo orbifolds.

Finally we mention that many new \Se structures can be constructed in higher
dimension by applying the join construction introduced in [BG1] to our new
examples. For example, 

\noindent{\sc Proposition} \linkpr.2: \tensl The joins $S^3\star L_{60}$ and
$S^3\star L_{256}^{(i)}$ are all smooth \Se 7-manifolds. $S^3\star L_{60}$ has
the rational cohomology type of $S^2\times \bigl((S^2\times S^3)\#(S^2\times
S^3)\bigr),$ while $S^3\star L_{256}^{(i)}$ have the rational cohomology type
of $S^3\star S^3\star S^3.$ \tenrm

\noindent{\sc Proof}: By Proposition 4.6 of [BG1] these joins are smooth
when the Fano indices of the links are one. But it follows from the proof of
Proposition \linkke.5 that $\gro^{-1}(L_f)\approx \calo(|\bfw|-d)$ and 
in all three cases we have $|\bfw|-d=1.$ The rational cohomology types can
easily be seen from Theorem 5.22 of [BG1]. \hfill\za

Many other examples can be worked out along the lines of [BG1]. However, what
is perhaps a more interesting question is, for example, whether  $S^3\star
L_{256}^{(1)}, S^3\star L_{256}^{(2)}$ and $S^3\star S^3\star S^3$ have the
same integral cohomology type, and if they do, are they homeomorphic
(diffeomorphic)? We plan to study these types of questions in the future.

\bigskip
\medskip
\centerline{\bf Bibliography}
\medskip
\font\ninesl=cmsl9
\font\bsc=cmcsc10 at 10truept
\parskip=1.5truept
\baselineskip=11truept
\ninerm

\item{[Arn]} {\bsc V.I. Arnold}, {\ninesl Some remarks on symplectic monodromy
of Milnor fibrations}, in Floer Memorial Volume, Progress in Mathematics 133,
Birkh\"auser, 1995, H. Hofer, C.H. Taubes, A. Weinstein, E. Zehnder Eds.
\item{[Bai]} {\bsc W. L. Baily}, {\ninesl On the imbedding of V-manifolds in
projective space}, Amer. J. Math. 79 (1957), 403-430.
\item{[BG1]} {\bsc C. P. Boyer and  K. Galicki}, {\ninesl On Sasakian-Einstein
Geometry}, Int. J. Math. 11 (2000), 873-909.
\item{[BG2]} {\bsc C. P. Boyer and  K. Galicki}, {\ninesl 
3-Sasakian manifolds}. {\it Surveys in differential geometry: 
essays on Einstein manifolds}, 123--184, Surv. Differ. Geom.,
VI, C. LeBrun and M. Wang, Eds., Int. Press, Boston, MA, 1999.
\item{[BGN1]} {\bsc C. P. Boyer, K. Galicki, and M. Nakamaye}, {\ninesl On the
Geometry of Sasakian-Einstein 5-Manifolds}, submitted
for publication; math.DG/0012041.
\item{[BGN2]} {\bsc C. P. Boyer, K. Galicki, and M. Nakamaye}, {\ninesl
Sasakian-Einstein Structures on $\scriptstyle{9\#(S^2\times S^3)}$}, submitted
for publication; math.DG/0102181.
\item{[BM]} {\bsc S. Bando and T. Mabuchi}, {\ninesl Uniqueness of Einstein
K\"ahler Metrics Modulo Connected Group Actions}, Adv. Stud. Pure Math. 10
(1987), 11-40.
\item{[B\"oh]} {\bsc C. B\"ohm}, {\ninesl Inhomogeneous Einstein metrics on
low dimensional spheres and other low dimensional spaces}, Invent. Math. 134
(1998), 145-176.
\item{[BR]} {\bsc M. Beltrametti and L. Robbiano} {\ninesl Introduction to the theory of 
weighted projective spaces}, Expo. Math. 4 (1986), 111-162.
\item{[Del]} {\bsc Ch. Delorme}, {\ninesl Espaces projectifs anisotropes},
Bull. Soc. Math. France 103 (1975), 203-223.
\item{[DK]} {\bsc J.-P. Demailly and J. Koll\'ar}, {\ninesl Semi-continuity of
complex singularity exponents and K\"ahler-Einstein metrics on Fano
orbifolds}, preprint AG/9910118.
\item{[Dim]} {\bsc A. Dimca}, {\ninesl Singularities and Topology of
Hypersurfaces}, Springer-Verlag, New York, 1992.
\item{[DT]} {\bsc W. Ding and G. Tian}, {\ninesl K\"ahler-Einstein metrics and
the generalized Futaki invariants}, Invent. Math. 110 (1992), 315-335.
\item{[Dol]} {\bsc I. Dolgachev}, {\ninesl Weighted projective varieties}, in
Proceedings, Group Actions and Vector Fields, Vancouver (1981) LNM 956, 34-71. 
\item{[DuKa]} {\bsc A. Durfee and L. Kauffman}, {\ninesl Periodicity of
Branched Cyclic Covers}, Math. Ann. 218 (1975), 175-189.
\item{[Fle]} {\bsc A.R. Fletcher}, {\ninesl Working with weighted complete
intersections}, Preprint MPI/89-95.
\item{[FK]} {\bsc Th. Friedrich and I. Kath}, {\ninesl Einstein manifolds of
dimension five with small first eigenvalue of the Dirac operator}, J. Diff.
Geom. 29 (1989), 263-279.
\item{[Hae]} {\bsc A. Haefliger}, {\ninesl Groupoides d'holonomie et
classifiants}, Ast\'erisque 116 (1984), 70-97.
\item{[KMM]} {\bsc Y. Kawamata, K. Matsuda, and K. Matsuki}, {\ninesl
Introduction to the Minimal Model Problem}, Adv. Stud. Pure Math. 10 (1987),
283-360.
\item{[Kol]} {\bsc J. Koll\'ar}, {\ninesl Rational Curves on Algebraic
Varieties}, Springer-Verlag, New York, 1996.
\item{[Mil]} {\bsc J. Milnor}, {\ninesl Singular Points of Complex
Hypersurfaces}, Ann. of Math. Stud. 61, Princeton Univ. Press, 1968.
\item{[MO]} {\bsc J. Milnor and P. Orlik}, {\ninesl Isolated singularities
defined by weighted homogeneous polynomials}, Topology 9 (1970), 385-393. 
\item{[Mo]} {\bsc S. Mori}, {\ninesl On a generalization of complete
intersections}, J. Math. Kyoto Univ. 15-3 (1975), 619-646.
\item{[Mor]} {\bsc S. Moroianu}, {\ninesl Parallel and Killing spinors on
$\hbox{Spin}^c$-manifolds}, Commun. Math. Phys. 187 (1997), 417-427.
\item{[MS]} {\bsc D. McDuff and D. Salamon}, {\ninesl Introduction to Symplectic Topology},
Oxford Mathematical Monographs, Oxford Univ. Press, Oxford, 1995.
\item{[Na]} {\bsc A.M. Nadel}, {\ninesl Multiplier ideal sheaves
and existence of K\"ahler-Einstein metrics of positive scalar curvature}, Ann.
Math. 138 (1990), 549-596.
\item{[Or1]} {\bsc P. Orlik}, {\ninesl Weighted homogeneous polynomials and
fundamental groups}, Topology 9 (1970), 267-273. 
\item{[Or2]} {\bsc P. Orlik}, {\ninesl On the homology of weighted homogeneous
manifolds}, Proc. 2nd Conf. Transformations Groups I, LNM 298,
Springer-Verlag, (1972), 260-269.
\item{[Ran]} {\bsc R.C. Randell} {\ninesl The homology of generalized
Brieskorn manifolds}, Topology 14 (1975), 347-355.
\item{[Siu]} {\bsc Y.-T. Siu}, {\ninesl The existence of K\"ahler-Einstein
metrics on manifolds with positive anticanonical line bundle and a suitable
finite symmetry group}, Ann. Math. 127 (1988), 585-627.
\item{[Sm]} {\bsc S. Smale}, {\ninesl On the structure of 5-manifolds},
Ann. Math. 75 (1962), 38-46.
\item{[Tak]} {\bsc T. Takahashi}, {\ninesl Deformations of Sasakian structures
and its applications to the Brieskorn manifolds}, Tohoku Math. J. 30 (1978),
37-43.
\item{[Tan]} {\bsc S. Tanno}, {\ninesl On the isometry groups of Sasakian
manifolds}, J. Math. Soc. Japan 22 (1970), 579-590.
\item{[Ti1]} {\bsc G. Tian}, {\ninesl On K\"ahler-Einstein metrics on certain
K\"ahler manifolds with $C_1(M)>0$}, Invent. Math. 89 (1987), 225-246.
\item{[Ti2]} {\bsc G. Tian}, {\ninesl On Calabi's Conjecture for complex
surfaces with positive first Chern class}, Invent. Math. 101 (1990), 101-172.
\item{[Ti3]} {\bsc G. Tian}, {\ninesl K\"ahler-Einstein metrics with positive
scalar curvature}, Invent. Math. 137 (1997), 1-37.
\item{[TY]} {\bsc G. Tian and S.-T. Yau}, {\ninesl K\"ahler-Einstein
metrics on complex surfaces with $c_1>0$}, Comm. Math. Phys. 112 (1987),
175-203.
\item{[Wa]} {\bsc M. Wang}, {\ninesl Einstein Metrics from Symmetry and Bundle
Constructions}, to appear in Essays on Einstein Manifolds, Surveys in 
Differential
Geometry Vol V, International Press, 2000, C. LeBrun and M. Wang, Eds. 
\item{[WZ]} {\bsc M. Wang and W. Ziller}, {\ninesl Einstein metrics on
principal torus bundles}, J. Diff. Geom. 31 (1990), 215-248. 
\item{[YK]} {\bsc K. Yano and M. Kon}, {\ninesl
Structures on manifolds}, Series in Pure Mathematics 3, 
World Scientific Pub. Co., Singapore, 1984.

\medskip
\bigskip \line{ Department of Mathematics and Statistics
\hfil March 2000} \line{ University of New Mexico \hfil revised September 2001} 
\line{ Albuquerque, NM 87131 \hfil } \line{ email: cboyer@math.unm.edu,
galicki@math.unm.edu\hfil} \line{ web pages:
http://www.math.unm.edu/$\tilde{\phantom{o}}$cboyer, 
http://www.math.unm.edu/$\tilde{\phantom{o}}$galicki \hfil}

\bye